\documentclass[11pt]{article}
\renewcommand{\baselinestretch}{1.2}
 \usepackage{amsmath,amssymb,amscd}
\newcommand{\pf}[1]{\trivlist \item[\hskip\labelsep\it #1\ ]}
\newcommand{\varpf}[1]{\trivlist \item[\hskip\labelsep\sc #1:]}
\newcommand{\qedbox}{$\rlap{$\sqcap$}\sqcup$}
\newcommand{\qed}{\qquad \qedbox \endtrivlist}
\newcommand{\varqed}{\hfill \rule{0.6em}{0.6em} \endtrivlist}
\newenvironment{proof}{\pf{Proof}}{\qed}

\newenvironment{items}{
  \begin{enumerate} 
                    
  }{\end{enumerate}}

\newenvironment{keywords}{\noindent\small {\it Keywords\/}:}{\vskip 4pt}
\newenvironment{classification}{\noindent\small 2000 {\it Mathematics Subject
Classification\/}:}{\vskip 12pt}

\newcommand{\tensor}{\otimes}

\newcommand{\Tensor}{\hat{\otimes}}

\newcommand{\wTensor}{\check{\otimes}}

\newcommand{\cstar}{{C^\ast}}

\newcommand{\id}{{\mathrm{id}}}

\newcommand{\A}{{\mathfrak A}}
\newcommand{\B}{{\mathfrak B}}

\newcommand{\op}{{\mathrm{op}}}

\newcommand{\El}{{{\cal E}\!\ell}}

\newtheorem{theorem}{Theorem}[section]
\newtheorem{lemma}[theorem]{Lemma}
\newtheorem{corollary}[theorem]{Corollary}
\newtheorem{proposition}[theorem]{Proposition}
\newtheorem{df}[theorem]{Definition}

\title{The flip is often discontinuous}
\author{\it Volker Runde}
\date{}
\begin{document}
\maketitle
\begin{abstract}
Let $\A$ be a Banach algebra. The flip on $\A \tensor \A^\op$ is defined through $\A \tensor \A^\op \ni a \tensor b 
\mapsto b \tensor a$. If $\A$ is ultraprime, $\El(\A)$, the algebra of all elementary operators on $\A$, can be
algebraically identified with $\A \tensor \A^\op$, so that the flip is well defined on $\El(\A)$. We show that
the flip on $\El(\A)$ is discontinuous if $\A = {\cal K}(E)$ for a reflexive Banach space $E$ with the approximation
property.
\end{abstract}
\section{The problem and the result}
A linear operator $T$ on an algebra $\A$ is called {\it elementary\/} if there are $a_1, b_1, \ldots, a_n, b_n \in \A$ such
that
\begin{equation} \tag{\mbox{$\ast$}} \label{op}
  Tx = \sum_{j=1}^n a_j x b_j \qquad (x \in \A).
\end{equation}
The set of all elementary operators on $\A$ is denoted by $\El(\A)$; equipped with the composition of operators as 
multiplication, $\El(\A)$ is an algebra. Let $\A^\op$ denote the opposite algebra of $\A$, i.e.\ the algebra we obtain by 
reversing the multiplication on $\A$. The canonical map that assigns to an element 
$\sum_{j=1}^n a_j \tensor b_j \in \A \tensor \A^\op$ the operator $T$ defined in (\ref{op}) is an algebra homomorphism 
from $\A \tensor \A^\op$ onto $\El(\A)$.
\par
Let $\A$ be a Banach algebra. Then $\El(\A)$ consists of bounded operators, and thus is a normed algebra in a natural
way. If $\A \tensor \A^\op$ is equipped with the projective norm, the canonical map from $\A \tensor \A^\op$ onto
$\El(\A)$ is continuous. If $\A$ is {\it ultraprime\/}, i.e.\ if there is $C > 0$ such that
\[
  \sup\{ \| a x b \| : x \in \A, \, \| x \| \leq 1 \} \geq C \| a \| \| b \| \qquad (a,b \in \A),
\]
the canonical map is even an (algebraic) isomorphism (\cite[Theorem 5.1]{Mat}). Examples of ultraprime Banach algebras
are all prime $\cstar$-algebras (\cite{Mat0}), as well as, for an arbitrary Banach space $E$, all closed subalgebras of 
${\cal B}(E)$ containing the finite rank operators and all closed ideals of factors (\cite[p.\ 305]{Mat}).
\par
The flip on $\A \tensor \A^\op$ is defined through $\A \tensor \A^\op \ni a \tensor b \mapsto b \tensor a$. It is
obviously continuous with respect to the projective norm. If $\A$ is ultraprime, the normed algebras $\A \tensor \A^\op$
and $\El(\A)$ are algebraically isomorphic, so that the flip is well defined on $\El(\A)$. It has been an open problem 
whether the flip on $\El(\A)$ is continuous (with respect to the operator norm). 
\par
It is the purpose of this note to give a negative answer to this question for certain ultraprime Banach algebras:
\begin{theorem} \label{thm}
Let $E$ be a reflexive Banach space with the approximation property. Then the flip on $\El({\cal K}(E))$ is
continuous if and only if $\dim E < \infty$.
\end{theorem}
\par
Let $E$ be a Banach space as in Theorem \ref{thm}, and let $\A$ be any closed subalgebra of ${\cal B}(E)$ containing
${\cal K}(E)$ (so that $\A$ is ultraprime). There is a canonical embedding of $\El({\cal K}(E))$ into $\El(\A)$.
Since $E^\ast$ has the bounded approximation property, ${\cal K}(E)$ has a bounded approximate identity
(\cite[Theorem 3.3]{GW}). Consequently, the operator norm on $\El({\cal K}(E))$ and the norm it inherits as a subalgebra
of $\El(\A)$ are equivalent. Since $\El({\cal K}(E))$ viewed as a subalgebra of $\El(\A)$ is invariant under the flip,
we get the following extension of Theorem \ref{thm}:
\begin{corollary}
Let $E$ be a reflexive Banach space with the approximation property, and let $\A$ be a closed subalgebra of ${\cal B}(E)$ 
containing ${\cal K}(E)$. Then the flip on $\El(\A)$ is continuous if and only if $\dim E < \infty$.
\end{corollary} 
\section{The proof}
If $\A$ is a normed algebra, its second dual can be equipped with two natural products extending the product on $\A$: the
first and the second Arens product (see \cite[1.4]{Pal} for details).
\par
The following proposition is essentially the discussion on \cite[p.\ 49]{Gro}. For any Banach space $E$, let
${\cal I}(E)$ denote the integral operators and ${\cal N}(E)$ the nuclear operators on $E$; for the composition of
operators we write $\circ$.
\begin{proposition} \label{prop}
Let $E$ be a Banach space such that $E^\ast$ has the bounded approximation property and ${\cal I}(E^\ast) = {\cal N}(E^\ast)$.
Then:
\begin{items}
\item ${\cal K}(E)^{\ast\ast}$ and ${\cal B}(E^{\ast\ast})$ are canonically isomorphic as Banach spaces.
\item The first Arens product $\circ_1$ on ${\cal K}(E)^{\ast\ast}$ is given by
\[
  S \circ_1 T = S \circ T \qquad (S, T \in {\cal B}(E^{\ast\ast})).
\]
\item The second Arens product $\circ_2$ on ${\cal K}(E)^{\ast\ast}$ is given by
\[
  S \circ_2 T = (j^\ast \circ S^{\ast\ast} \circ i^{\ast\ast}) \circ T \qquad (S,T \in {\cal B}(E^{\ast\ast})),
\]
where $i  \!: E \to E^{\ast\ast}$ and $j \!: E^\ast \to E^{\ast\ast\ast}$ are the canonical
embeddings.
\end{items}
\end{proposition}
\par
Note that ${\cal I}(E^\ast) = {\cal N}(E^\ast)$ is always satisfied if $E^\ast$ has also the Radon--Nikod\'ym property
(\cite[16.6, Theorem]{DF}). It is an immediate consequence of Proposition \ref{prop}(i), that ${\cal K}(E)^{\ast\ast}$
equipped with the first Arens product has an identity. For the second Arens product, we have the following:
\begin{corollary} \label{cor}
Let $E$ be a Banach space such that $E^\ast$ has the bounded approximation property and ${\cal I}(E^\ast) = {\cal N}(E^\ast)$.
Then ${\cal K}(E)^{\ast\ast}$ equipped with the second Arens product has an identity if and only if $E$ is reflexive.
\end{corollary}
\begin{proof}
If $E$ is reflexive, ${\cal K}(E)$ is Arens regular (\cite[Corollary, p.\ 103]{Pal}), so that $\circ_1$ and $\circ_2$
coincide.
\par
For the converse, first note that, by Proposition \ref{prop}(iii), $\id_{E^{\ast\ast}}$ is a left identity for 
$({\cal K}(E)^{\ast\ast}, \circ_2)$. Hence, if $({\cal K}(E)^{\ast\ast}, \circ_2)$ has an identity, it must be 
$\id_{E^{\ast\ast}}$. From Proposition \ref{prop}(iii), we thus obtain
\[
  S = j^\ast \circ S^{\ast\ast} \circ i^{\ast\ast} = (i^\ast \circ S^\ast \circ j)^\ast \qquad (S \in {\cal B}(E^{\ast\ast})).
\]
This means that every bounded, linear operator on $E^{\ast\ast}$ is the adjoint of a bounded, linear operator on $E^\ast$,
which is possible only if $E^\ast$ is reflexive. It follows that $E$ is reflexive.
\end{proof}
\par
The following lemma is certainly well known, but we could find it nowhere in the literature; it follows immediately from
the separate $w^\ast$-continuity properties of the first and the second Arens product:
\begin{lemma} \label{lemma}
Let $\A$ and $\B$ be Banach algebras, and let $\theta \!: \A \to \B$ be a continuous anti-homomorphism. Then, if
$\A^{\ast\ast}$ is equipped with the first Arens product and $\B^{\ast\ast}$ is equipped with the second Arens product,
$\theta^{\ast\ast} \!: \A^{\ast\ast} \to \B^{\ast\ast}$ is also an anti-homomorphism.
\end{lemma}
\par
Obviously, if $\theta$ is an anti-isomorphism, then so is $\theta^{\ast\ast}$. As an immediate consequence, we obtain:
\begin{corollary} \label{cor2}
Let $\A$ be a Banach algebra which is topologically anti-isomorphic to itself. Then $\A^{\ast\ast}$ has an identity with 
respect to the first Arens product if and only if it has an identity with respect to the second Arens product.
\end{corollary}
\par
This can be used to rule out the existence of certain anti-isomorphisms:
\begin{lemma} \label{lemma2}
Let $E$ be a reflexive, infinite-dimensional Banach space with the approximation property. Then ${\cal K}({\cal K}(E))$ is 
not topologically anti-isomorphic to itself.
\end{lemma}
\begin{proof}
First note that, since $E$ is reflexive, it does in fact have the metric approximation property (\cite[16.4, Corollary 4]{DF}).
Since $E$ is reflexive, the same is true for $E^\ast$ (\cite[5.7, Corollary]{DF}). Since 
\[
  {\cal K}(E)^\ast \cong {\cal N}(E^\ast) \cong E^\ast \Tensor E,
\]
it is easily checked that ${\cal K}(E)^\ast$ has also the metric approximation property. Moreover, ${\cal N}(E^\ast)$
has the Radon--Nikod\'ym property (\cite[p.\ 219]{DU}), so that ${\cal I}({\cal K}(E)^\ast) = {\cal N}({\cal K}(E)^\ast)$, 
i.e.\ the hypotheses of Proposition \ref{prop} are satisfied. Assume towards a contradiction that ${\cal K}({\cal K}(E))$
is topologically anti-isomorphic to itself. Since ${\cal K}({\cal K}(E))^{\ast\ast}$ equipped with the first Arens
product has an identity, Corollary \ref{cor2} implies that the same is true for ${\cal K}({\cal K}(E))^{\ast\ast}$ equipped 
with the second Arens product. This, however, means that ${\cal K}(E)$ is reflexive by Corollary \ref{cor}. Since
${\cal K}(E)^{\ast\ast} \cong {\cal B}(E)$, again by Proposition \ref{prop}, it follows that every bounded, linear
operator on $E$ must be compact, i.e.\ $\dim E < \infty$.
\end{proof}
\par
By \cite[Theorem 33.3(i)]{BD}, $\El({\cal K}(E))$ is contained in ${\cal K}({\cal K}(E))$ for any Banach space $E$. As well 
shall see, for the spaces $E$ we are interested in, $\El({\cal K}(E))$ is even dense in ${\cal K}({\cal K}(E))$:
\begin{proposition} \label{prop2}
Let $E$ be a reflexive Banach space with the approximation property. Then $\El({\cal K}(E))$ is norm dense in
${\cal K}({\cal K}(E))$.
\end{proposition}
\par
As preparation for the proof, we need two lemmas, which are certainly known, but for which we could not find an explicit
reference. For convenience, we include proofs.
\par
For any normed space $E$, we write ${\cal F}(E)$ for the bounded finite rank operators on $E$.
\begin{lemma} \label{l1}
Let $E$ be a normed space, and let $T \in {\cal F}(E)$. Then there is a projection $P \in {\cal F}(E)$ such that
$PT = TP = T$.
\end{lemma}
\begin{proof}
Let $x_1, \ldots, x_n \in E$ and $\phi_1, \ldots, \phi_n \in E^\ast$ be such that
\[
  T = \sum_{j=1}^n x_j \tensor \phi_j.
\]
Choose finite-dimensional subspaces $X$ of $E$ and $Y$ of $E^\ast$ with $x_1, \ldots, x_n \in X$ and $\phi_1, \ldots, \phi_n
\in Y$ such that $(X,Y)$ forms a dual systems, i.e.\ the bilinear form $X \times Y \ni (x, \phi) \mapsto \langle x , \phi 
\rangle$ is non-degenerate; note that necessarily $\dim X = \dim Y$. Choose bases $x^\prime_1, \ldots x^\prime_m$ of
$X$ and $\phi^\prime_1, \ldots, \phi_m^\prime$ of $Y$ such that
\[
  \langle x^\prime_j, \phi^\prime_k \rangle = \delta_{j,k} \qquad (j,k =1, \ldots, m).
\]
Define
\[
  P := \sum_{j=1}^m x_j^\prime \tensor \phi^\prime_j.
\]
It is easily seen that $PT = TP = T$.
\end{proof}
\begin{lemma} \label{l2}
Let $E$ be a normed space, and let $F$ be a finite-dimensional subspace of ${\cal F}(E)$. Then there is a projection $p
\in {\cal F}(E)$ such that $F \subset p {\cal F}(E) p$.
\end{lemma}
\begin{proof}
Since ${\cal F}(E) \cong E \tensor E^\ast$, we may find finite-dimensional subspaces $X$ of $E$ and $Y$ of
$E^\ast$ such that $F \subset X \tensor Y$. Making $X$ and $Y$ larger, if necessary, we may suppose that $(X,Y)$
forms a dual system. Defining $p$ analoguously to $P$ in the proof of Lemma \ref{l1}, we obtain the desired projection.
\end{proof}
\pf{Proof of Proposition\/ {\rm \ref{prop2}}}
Let $T \in {\cal K}({\cal K}(E))$. As have already seen in the proof of Lemma \ref{lemma2}, ${\cal K}(E)^\ast$ has the
(metric) approximation property, so that ${\cal K}({\cal K}(E)) \cong {\cal K}(E) \wTensor {\cal K}(E)^\ast$. Since
${\cal F}(E) \tensor {\cal K}(E)^\ast$ is norm dense in ${\cal K}(E) \wTensor {\cal K}(E)^\ast$, we may suppose that
$T \in {\cal F}(E) \tensor {\cal K}(E)^\ast \cong {\cal F}({\cal F}(E))$. Let $P \in {\cal F}({\cal F}(E))$ as specified
in Lemma \ref{l1}. Furthermore, since $\dim P{\cal F}(E) < \infty$, Lemma \ref{l2}, yields a projection $p \in
{\cal F}(E)$ such that $P{\cal F}(E) \subset p {\cal F}(E) p = p {\cal K}(E) p =: \A$. Then $\A$ is a subalgebra
of ${\cal K}(E)$ such that $T \A \subset \A$. Obviously, $\A \cong {\cal B}(pE)$. Since ${\cal B}(pE)$ is ultraprime
and finite-dimensional, it follows from \cite[Theorem 5.1]{Mat} and the elementary fact that $\A \tensor \A^\op$
and ${\cal B}(\A)$ have the same finite dimension that $\El(\A) = {\cal B}(\A)$. Hence, there are $a_1, b_1, \ldots,
a_n , b_n \in \A$ such that
\[
   Tx = \sum_{j=1}^n a_j x b_j \qquad (x \in \A).
\]
From the choice of $p$, it follows that $T(pxp) = Tx$ for all $x \in {\cal K}(E)$, so that 
\[
  Tx = T(pxp) =  \sum_{j=1}^n a_jpxp b_j = \sum_{j=1}^n a_j x b_j \qquad (x \in {\cal K}(E)), 
\]
i.e.\ $T \in \El({\cal K}(E))$.
\qed
\par
With all these preparations made, a proof of Theorem \ref{thm} is now a matter of a few lines:
\pf{Proof of Theorem\/ {\rm \ref{thm}}}
First, observe that the flip is an anti-automorphism of $\El({\cal K}(E))$ which is its own inverse. Suppose 
that the flip is continuous on $\El({\cal K}(E))$. By Proposition \ref{prop2}, it then extends to a topological
anti-automorphism of ${\cal K}({\cal K}(E))$. By Lemma \ref{lemma2}, this is possible only if $E$ is finite-dimensional.
\qed 
\vfill
\renewcommand{\baselinestretch}{1.2}
\begin{tabbing} 
{\it Address\/}: \= Department of Mathematical Sciences \\
                 \> University of Alberta \\
                 \> Edmonton, Alberta \\
                 \> Canada T6G 2G1 \\ \medskip
                 {\it E-mail\/}:  \> {\tt runde@math.ualberta.ca}
\end{tabbing}


\begin{thebibliography}{0}
%
\begin{small}
%
\bibitem{BD} {\sc F.\ F.\ Bonsall} and {\sc J.\ Duncan}, {\it Complete Normed Algebras\/}. Springer Verlag, 1973.
%
\bibitem{DF} {\sc A.\ Defant} and {\sc K.\ Floret}, {\it Tensor Norms and Operator Ideals\/}. North-Holland, 1993.
%
\bibitem{DU} {\sc J.\ Diestel} and {\sc J.\ J.\ Uhl, Jr.}, {\it Vector Measures\/}. American Mathematical Society, 1977.
%
\bibitem{GW} {\sc N.\ Gr{\oe}nb{\ae}k} and {\sc G.\ A.\ Willis}, Approximate identities in Banach algebras of compact
operators. {\it Canad.\ Math.\ Bull.\/}\ {\bf 36\/} (1993), 45--53.
%
\bibitem{Gro} {\sc M.\ Grosser}, Arens semiregular Banach algebras. {\it Monatsh.\ Math.\/}\ {\bf 98\/} (1984),
41--52.
%
\bibitem{Mat0} {\sc M.\ Mathieu}, Elementary operators on prime $\cstar$-algebras, I. {\it Math.\ Ann.\/}\ {\bf 284\/} (1989),
223--244.
%
\bibitem{Mat} {\sc M.\ Mathieu}, Rings of quotients of ultraprime Banach algebras, with applications to elementary operators.
In: {\sc R.\ J.\ Loy} (ed.), {\it Conference on Automatic Continuity and Banach Algebras\/}. Australian National University,
1989, pp.\ 297--317.
%
\bibitem{Pal} {\sc T.\ W.\ Palmer\/}, {\it Banach Algebras and the General Theory of $^\ast$-Algebras\/}, I. Cambridge 
University Press, 1994.     
%
\end{small}
%
\end{thebibliography}
\end{document}